\documentclass[10pt]{amsart}
\usepackage{amsfonts,amssymb,amscd,amsmath,enumerate,verbatim,calc}
\usepackage{color,soul}
\usepackage{mathtools}
\usepackage{graphicx}
\everymath{\displaystyle}
\usepackage{fancyhdr}
\usepackage{tikz}
\usepackage{url}
\usepackage{amsmath}
\usetikzlibrary{shapes.geometric, arrows}
\tikzstyle{arrow} = [thick,->,>=stealth]
\tikzstyle{process} = [rectangle, minimum width=4cm, minimum height=2cm, text centered, text width=2.5cm, draw=green, fill=green!10]
\theoremstyle{plain}

\theoremstyle{definition}

\begin{document}
\title{SOME SEPARATION AXIOMS IN TOPOLOGICAL SPACES}
\author{Neeraj kumar Tomar}
\email{neer8393@gmail.com}
\address{Department of Applied Mathematics, Gautam Buddha University, Greater Noida, Uttar Pradesh 201312, India}

\author{M. C. Sharma}
\email{sharmamc2@gmail.com}
\address{Professor Department of Mathematics, N. R. E. C. College Khurja, Uttar Pradesh 203131, India}

\author{Amit Ujlayan}
\email{amitujlayan@gbu.ac.in}
\address{Department of Applied Mathematics, Gautam Buddha University, Greater Noida, Uttar Pradesh 201312, India}

\author{Fahed Zulfeqarr}
\email{fahed@gbu.ac.in}
\address{Department of Applied Mathematics, Gautam Buddha University, Greater Noida, Uttar Pradesh 201312, India}


\date{\today}
	
\begin{abstract}
In this paper, we introduced the concepts of new separation axioms called \( SC^* \)-separation axioms and \( H^* \)-separation axioms by using \( SC^* \) and \( H^* \)-open sets in topological spaces. The \( SC^* \)-separation axioms include \( SC^* \)-\( C_0 \), \( SC^* \)-\( C_1 \), weakly \( SC^* \)-\( C_0 \), and weakly \( SC^* \)-\( C_1 \) spaces, while the \( H^* \)-separation axioms include \( H^* \)-\( T_{\frac{1}{2}} \), \( H^* \)-\( T_b \), and \( H^* \)-\( T_d \) spaces. Also, we obtained several properties of such spaces.
\end{abstract}
\maketitle
\section{Introduction}
In this paper, we introduced the concepts of new separation axioms called \( SC^* \)- and \( H^* \)-separation axioms, i.e., \( SC^* \)-\( C_0 \), \( SC^* \)-\( C_1 \), weakly \( SC^* \)-\( C_0 \), weakly \( SC^* \)-\( C_1 \), and \( H^* \)-\( T_{\frac{1}{2}} \), \( H^* \)-\( T_b \), \( H^* \)-\( T_d \), by using \( SC^* \)- and \( H^* \)-open sets. These developments are inspired by the works of A. Chandrakala and K. Bala Deepa Arasi~\cite{chandrakala2022scclosed}, and Nidhi Sharma et al.~\cite{sharma2015h} in the context of topological spaces. We also obtained several properties of such spaces.
In 2024, K. Suthi Keerthana et al.~\cite{keerthana2024abgaclosure} introduced a new class of separation axioms based on \( \alpha b^*g\alpha \)-closed sets, namely the \( cT_{\alpha b^*g\alpha} \)-space, \( T^*_{\alpha b^*g\alpha} \)-space, and \( *_{g\alpha} T_{1/2}^{***} \)-space, using \( \alpha b^*g\alpha \)-open and \( b^*g\alpha \)-open sets. They explored several properties of these generalized spaces.

\section{Prerequisites}
Throughout this paper, by a space \( X \), we mean a topological space \( (X, \tau) \). If \( A \) is any subset of a space \( X \), then \( \text{cl}(A) \), \( \text{int}(A) \), and \( C(A) \) denote the closure, the interior, and the complement of \( A \), respectively. The family of all semi-open (resp. pre-open, \( \alpha \)-open, \( c^* \)-open) subsets of a space \( X \) is denoted by \( SO(X) \) (resp. \( PO(X) \), \( \alpha O(X) \), \( c^*O(X) \)).
Maheswari and Tapi~\cite{maheswari1978feeblyT1} called a subset \( B \) of a space \( X \) \emph{feebly open} if there exists an open set \( G \) such that \( G \subseteq B \subseteq \text{scl}(G) \). Later, Janković and Reilly~\cite{jankovic1985semiseparation} observed that feebly open sets are precisely the \( \alpha \)-open sets.

Let us recall the following definitions:

\subsection{Definition:}

A subset \( A \) of a topological space \( X \) is said to be:
\begin{enumerate}
    \item \textbf{semi-open}\cite{levine1963semi} if \( A \subseteq \mathrm{cl}(\mathrm{int}(A)) \).
    \item \textbf{\( \alpha \)-open} \cite{njȧstad1965some} if \( A \subseteq \mathrm{int}(\mathrm{cl}(\mathrm{int}(A))) \)
    \item \textbf{pre-open}~\cite{mashhour1982precontinuous} if \( A \subseteq \mathrm{int}(\mathrm{cl}(A)) \).
    \item \textbf{\( c^* \)-open}~\cite{malathi2017pregencstar} if 
    \[
    \mathrm{int}(\mathrm{cl}(A)) \subseteq A \subseteq \mathrm{cl}(\mathrm{int}(A)).
    \]
    \item \textbf{\( SC^* \)-closed}~\cite{chandrakala2022scclosed} if \( \mathrm{scl}(A) \subseteq U \) whenever \( A \subseteq U \) and \( U \) is \( c^* \)-open in \( X \).
\end{enumerate}
The complement of a semi-open (resp. \( \alpha \)-open, pre-open, \( c^* \)-open) subset of a space \( X \) is called a semi-closed (resp. \( \alpha \)-closed, pre-closed, \( c^* \)-closed) set. The semi-closure (resp. pre-closure, \( \alpha \)-closure, \( c^* \)-closure) of a set \( A \) is denoted by \( \mathrm{scl}(A) \) (resp. \( \mathrm{pcl}(A) \), \( \alpha\mathrm{cl}(A) \), \( c^*\mathrm{cl}(A) \)).

The complement of an \( SC^* \)-closed set is called an \( SC^* \)-open set.

\subsection{Definition:}

A subset \( A \) of a topological space \( X \) is said to be:

\begin{enumerate}
    \item \textbf{\( \alpha^* \)-set} \cite{hatir1996decomposition} if \( \operatorname{int}(\operatorname{cl}(\operatorname{int}(A))) = \operatorname{int}(A) \).
    \item \textbf{\( C \)-set} \cite{hatir1996decomposition} if \( A = U \cap V \), where \( U \) is an open set and \( V \) is an \( \alpha^* \)-set in \( X \).
     \item \( \alpha Cg \)-closed~\cite{devi2005stronglyacg} if \( \alpha\text{-cl}(A) \subseteq U \), whenever \( A \subseteq U \) and \( U \) is a \( C \)-set in \( X \).
     \item \textbf{\( w \)-closed} \cite{Sundaram2000} if \( \operatorname{cl}(A) \subseteq U \), whenever \( A \subseteq U \) and \( U \) is semi-open in \( X \).
    \item \textbf{\( h \)-closed} \cite{rodrigo2012international} if \( s\text{-}\operatorname{cl}(A) \subseteq U \), whenever \( A \subseteq U \) and \( U \) is \( w \)-open in \( X \).
    \item \textbf{\( hCg \)-closed} \cite{rodrigo2012international} if \( h\text{-}\operatorname{cl}(A) \subseteq U \), whenever \( A \subseteq U \) and \( U \) is a \( C \)-set in \( X \).
    \item \( H^* \)-closed~\cite{sharma2015h} if \( h\text{-cl}(A) \subseteq U \), whenever \( A \subseteq U \) and \( U \) is \( hC_g \)-open in \( X \).
    \item \( H^*g \)-closed~\cite{sharma2015h} if \( H^*\text{-cl}(A) \subseteq U \), whenever \( A \subsetneq U \) and \( U \) is open in \( X \).
    \item \( gH^* \)-closed~\cite{sharma2015h} if \( H^*\text{-cl}(A) \subseteq U \), whenever \( A \subseteq U \) and \( U \) is \( H^* \)-open in \( X \).  
\end{enumerate}

\noindent
The complement of an \( \alpha Cg \)-closed (resp. \( w \)-closed, \( h \)-closed, \( hCg \)-closed, \( H^* \)-closed, \( H^*g \)-closed, \( gH^* \)-closed) set is said to be an \( \alpha Cg \)-open (resp. \( w \)-open, \( h \)-open, \( hCg \)-open, \( H^* \)-open, \( H^*g \)-open, \( gH^* \)-open) set. The intersection of all \( H^* \)-closed subsets of \( X \) containing \( A \) is called the \( H^* \)-closure of \( A \), and is denoted by \( H^*\text{-cl}(A) \). The union of all \( H^* \)-open subsets of \( X \) that are contained in \( A \) is called the \( H^* \)-interior of \( A \), and is denoted by \( H^*\text{-int}(A) \).\\
The family of all \( H^* \)-open (resp. \( H^* \)-closed) sets in a topological space \( X \) is denoted by \( H^*O(X) \) (resp. \( H^*C(X) \)).

\subsection{Definition}

A topological space \( X \) is called:

\begin{enumerate}
    \item \( C_0 \) (resp. semi-\( C_0 \)) if, for \( x, y \in X \), \( x \neq y \), there exists \( G \in \tau \) (resp. \( G \in SO(X) \)) such that \( \mathrm{cl}(G) \) (resp. \( \mathrm{scl}(G) \)) contains only one of \( x \) or \( y \), but not both.

    \item \( C_1 \) (resp. semi-\( C_1 \)) if, for \( x, y \in X \), \( x \neq y \), there exist \( G, H \in \tau \) (resp. \( G, H \in SO(X) \)) such that 
    \[
    x \in \mathrm{cl}(G),\ y \in \mathrm{cl}(H),\ x \notin \mathrm{cl}(H),\ y \notin \mathrm{cl}(G)
    \]
    (resp. replace \( \mathrm{cl} \) with \( \mathrm{scl} \)).

    \item \( w\text{-}C_0 \)~\cite{dimaio1985weakerR0} if \( \bigcap_{x \in X} \ker(x) = \emptyset \), where \( \ker(x) = \bigcap \{ G : x \in G \in \tau \} \).

    \item weakly semi-\( C_0 \) if, \( \bigcap_{x \in X} \mathrm{sker}(x) = \emptyset \), where \( \mathrm{sker}(x) = \bigcap \{ G : x \in G \in SO(X) \} \).

    \item \( R_0 \)~\cite{davis1961indexed} if \( \mathrm{cl}(\{x\}) \subseteq G \) whenever \( x \in G \in \tau \).

    \item semi-\( R_0 \)~\cite{maheswari1975r0spaces} if, for \( x \in G \in SO(X) \), \( \mathrm{scl}(\{x\}) \subseteq G \).

    \item weakly \( R_0 \)~\cite{dimaio1985weakerR0} if, \( \bigcap_{x \in X} \mathrm{cl}(\{x\}) = \emptyset \).

    \item weakly semi-\( R_0 \)~\cite{arya1990snormal} if, \( \bigcap_{x \in X} \mathrm{scl}(\{x\}) = \emptyset \).

    \item weakly pre-\( R_0 \) if, \( \bigcap_{x \in X} \mathrm{pcl}(\{x\}) = \emptyset \).

    \item weakly pre-\( C_0 \) if, \( \bigcap_{x \in X} \mathrm{pker}(x) = \emptyset \), where \( \mathrm{pker}(x) = \bigcap \{ G : x \in G \in PO(X) \} \).

    \item \( \alpha \)-space~\cite{dontchev1994superconnected} if every \( \alpha \)-open set in it is open.
\end{enumerate}
\noindent

Maheshwari and Prasad\cite{maheswari1975newseparation} introduced semi-$T_i$ $(i = 0, 1, 2)$ axiom, which is weaker than  $T_i$ $(i = 0, 1, 2)$ axiom.
\\
We use the following sets and classes for counter examples.
\\
Let $X$ = $\{a, b, c, d\}$,  \hspace{1.6cm} $Y$ = $\{a, b, c\}$, \hspace{1.5cm} $Z$ = $\{a, b, c, d, e,f\}$
\\
Let $\tau_1$ = $\{\phi,X, \{b\}, \{a,b\}, \{b,c\}, \{a,b,  c\}\}$, \hspace{1.3cm} $\sigma_1$ = $\{\phi, X, \{a\}, \{b\}, \{a,b\}\}$ 
\vspace{.2cm}
\\
$\eta_1$ = $\{\phi,Z, \{a,c,e\}, \{b,d.f\}\}$, \hspace{3.3cm} $\sigma_2$ = $\{\phi,Y, \{a\}, \{b\}, \{a,b\}, \{b,  c\}\}$
\vspace{.2cm}
\\
$\tau_2$ = $\{\phi,X, \{a\}, \{b\}, \{c\}, \{a,b\}, \{a,c\}, \{b,c\}, \{a,b,  c\}\}$
\vspace{.2cm}
\\
$\sigma$ = $\{\phi, Y, \{a\}, \{b\}, \{a,b\}\}$.
\subsection{Remark}
    $(X, \sigma)$ is semi-$C_0$ but not $C_0$, $(X,\tau_2)$ is semi-$C_1$, $C_0$ but not $C_1$. $(Y, \sigma_2)$ is semi-$T_0$ but not semi-$C_0$. $(Z, \eta_1)$ is an $\alpha$-space but not an $\alpha$-$C_0$.

\subsection{Theorem}
\begin{enumerate}
    \item Evey $C_1$ (semi-$C_1$) space is a $C_0$ (semi-$C_0$).
\item  Evey $C_0$ ($C_1$) space is a semi-$C_0$ (semi-$C_1$).
\item  Evey $R_0$ space is a weakly-$R_0$.
\item  Evey wealky $R_0$ space is a weakly semi-$R_0$.
\item  Evey semi $C_0$ (semi-$C_1)$ space is a semi-$T_0$(semi-$T_1$).
\end{enumerate}
\begin{proof}
     Omitted.
 \end{proof}
 \section{SC\(^*\)-\(C_0\), SC\(^*\)-\(C_1\), Weakly SC\(^*\)-\(C_0\), and Weakly SC\(^*\)-\(R_0\) Spaces}
\subsection{Definition}  A topological space \( X \) is called:
\begin{enumerate}
    \item SC$^*$-C$_0$ if, for $x, y \in X$, $x \ne y$, there exists $G \in SC^*(X)$ such that $SC^*\text{cl}(G)$ contains only one of $x$ and $y$ but not the other.

    \item SC$^*$-C$_1$ if, for $x, y \in X$, $x \ne y$, there exist $G, H \in SC^*(X)$ such that 
    $x \in SC^*\text{cl}(G)$, $y \in SC^*\text{cl}(H)$ but 
    $x \notin SC^*\text{cl}(H)$ and $y \notin SC^*\text{cl}(G)$.

    \item weakly SC$^*$-C$_0$ if, $\bigcap\limits_{x \in X} SC^*\text{ker}(x) = \emptyset$, where $SC^*\text{ker}(x) = \bigcap \{ G : x \in G \in SC^*(X) \}$.

    \item weakly SC$^*$-R$_0$ if, $\bigcap\limits_{x \in X} SC^*\text{cl}(\{x\}) = \emptyset$.
\end{enumerate}

\subsection{Theorem} 
\begin{enumerate}
    \item Every $SC^*$-$C_1$ space is $SC^*$-$C_0$.
\item Every $SC^*$-$C_0$($SC^*$-$C_1$) space is semi-$C_0$(semi-$C_1$).
\item Every weakly $SC^*$-$R_0$ space is weakly semi-$R_0$ and weakly pre-$R_0$.
\item Every $w$-$C_0$ space is weakly $SC^*$-$C_0$.
\item Every weakly $SC^*$-$C_0$ space is weakly semi-$C_0$ and weakly pre-$C_0$.
\item Every $SC^*$-$C_0$ ($SC^*$-$C_1$) space is semi-$T_0$ (semi-$T_1$).
\item Every weakly $R_0$ space is weakly $SC^*$-$R_0$.
\item weakly $SC^*$-$R_0$ness and weakly $SC^*$-$C_0$ness are independent notions.
\end{enumerate}

\subsection{Remark}
$(X, \tau_2)$ is $SC^*$-$C_0$ but not $SC^*$-$C_1$. $(X, \sigma_1)$ is semi-$C_0$, semi-$C_1$, but neither $SC^*$-$C_1$ nor $SC^*$-$C_0$. $(Y,\sigma)$ is weakly semi-$R_0$ but not weakly $SC^*$-$R_0$. $(X, \tau_2)$ is weakly $SC^*$-$C_0$ but not weakly $SC^*$-$R_0$. $(X, \tau_1)$ is weakly $SC^*$-$R_0$ but not weakly $SC^*$-$C_0$.

\subsection{Theorem}
A topological space $X$ is weakly $SC^*$-$R_0$ if and only if $SC^*$ker$(x) \neq X$ for each $x \in X$.
\begin{proof}
    \textbf{Necessity:} If there is some $x_0 \in  X$ with $SC^*ker(x_0)$ = $X$, then $X$ is the only $SC^*$-open set containing $x_0$. This implies that $x_0 \in SC^*cl(\{x\})$ for every $x\in X$. Hence $\bigcap _{x\in X}$ $SC^*cl(\{x\})\neq \phi$, a contradiction.
\\
\textbf{Sufficiency:} If $X$ is not weakly $SC^*$-$R_0$, then choose some $x_0\in$ $\bigcap _{x\in X}$ $SC^*ker(x).$ This implies that every $SC^*$- neighbourhood of $x_0$ contains every point of $X$. Hence  $SC^*ker(x_0)$ = $X$.
\end{proof}

\subsection{Theorem.}
A topological space $X$ is weakly $SC^*$-$C_0$ if and only if for each $x \in X$, there exists a proper $SC^*$-closed set containing $x_0$.

\begin{proof}
    \textbf{Necessity:} Suppose there is some $x_0 \in  X$ such that $X$ is the only $SC^*$-closed set containing $x_0$. Let $U$ be any proper $SC^*$-open subset of $X$ containing a point $x$. This implies that $C(U) \neq X$. Since $C(U)$ is S$C^*$-closed, we have $x_0 \in C(U)$. So $x_0 \in  U$. Thus $x_0\in$ $\bigcap _{x\in X}ker(x)$ for any point $x$ of $X$, is a contradiction.
\\
\textbf{Sufficiency:} If $X$ is not weakly $SC^*$-$C_0$, then choose $x_0\in $ $\bigcap _{x\in X}SC^*ker(x)$ so $x_0$ belongs to $SC^*ker(x)$ for any $x \in X$. This implies that $X$ is the only $SC^*$-open set, which contains the point $x_0$, is a contradiction.
\end{proof}

\subsection{Theorem.}
Every $SC^*$-$C_0$ ($SC^*$-$C_1$) space is weakly $SC^*$-$C_0$.
\begin{proof}
     If $x, y \in X$ such that $x \neq y$, where $X$ is a $SC^*$-$C_0$ space, then without loss of generality, we can assume that there exists $U \in SC^*(X)$ such that $x \in SC^*cl(U)$ but $y \notin SC^*cl(U)$. This implies that $U\neq \phi$. Hence we can choose $z$ in $U$.\\ Now,
 $SC^*ker(z) \cap SC^*ker(y) \subseteq U \cap (SC^*cl(U))\subseteq (SC^*cl(U)) \cap C(SC^*cl(U))$ = $\phi$.\\ Hence $\bigcap _{x\in X}ker(x)$ = $\phi$. Since every $SC^*$-$C_1$ space is also $SC^*$-$C_0$, it is also clear that every $SC^*$-$C_1$ space is weakly $SC^*$-$C_0$.
\end{proof}
\subsection{Remark.} The converse of the above theorem need not be true since $(Z, \eta_1)$ is weakly $SC^*$-$C_0$ but not $SC^*$-$C_0$. The space $(Y, \sigma)$ is both $SC^*$-$C_0$ and weakly $SC^*$-$C_0$ but not $SC^*$-$C_1$.

\subsection{Theorem} The property of being an $SC^*$-$C_0$ space is not hereditary.
\begin{proof}
     Consider the space $(Y, \sigma)$. Let $S$ = $\{a, c\}$ and $\sigma^*_1$ be the relative topology on $S$. It is easy to verify that $(Y, \sigma)$ is $SC^*$-$C_0$ but it’s subspace $(S, \sigma^*_1)$ is not $SC^*$-$C_0$.
 \end{proof}
\section{\( H^* \)-Separation Axioms}

\subsection{Remark.} Every $\alpha$-closed (resp. $\alpha$-open) set is $H^*$-closed (resp. $H^*$-open) set.\\
     For definitions stated above, we have the following diagram:

\vspace{.1cm}

    closed $\Longrightarrow$ $\alpha$-closed  $\Longrightarrow$ $g\alpha$-closed $\Longrightarrow$ $\alpha g$-closed\\

\hspace{2.2cm} $\Downarrow \hspace{2.2cm} \Downarrow \hspace{2.2cm} \Downarrow$\\

\hspace{1.1cm} $H^*$-closed $\Longrightarrow$ $gH^*$-closed $\Longrightarrow$ $H^*g$-closed
\vspace{.4cm}\\
However the converses of the above are not true as may be seen by the following examples:

\subsubsection{Example.}
Let $X$ = $\{a,b,c,d\}$ and $\tau$ = $\{\phi, X, \{a\}, \{b\}, \{a,b\}, \{a,b,c\}\}$. Then $A$ = $\{c\}$ is $\alpha$-closed set as well as $H^*$-closed set but not closed set in $X$.

\subsubsection{Example.}
Let $X$ = $\{a,b,c,d\}$ and $\tau$ = $\{\phi, X, \{a\}, \{b\}, \{a,b\}, \{a,b,c\}\}$. Then $A$ = $\{c\}$ is $\alpha$-closed set as well as $gH^*$-closed set but not closed set in $X$.

\subsection{Remark.}
\textbf{(i).} A subset $A$ of $X$ is $H^*g$-open in $X$ iff $F \subset \tau^{H^*}$-$int(A)$ whenever $F \subset A$ and $A$ is closed in $X$.\\     
\textbf{(ii).} A subset $A$ of $X$ is $gH^*$-closed (resp. $gH^*$-open) in $X$ iff $A$ is $g$-closed (resp. $g$-open) in $X$.

\section[\texorpdfstring{$H^*$-$T_1$, $H^*$-$T_{\frac{1}{2}}$, $H^*$-$T_b$ and $H^*$-$T_d$-Spaces}{H*-T1, H*-T1/2, H*-Tb and H*-Td-Spaces}]{\boldmath $H^*$-$T_1$, $H^*$-$T_{\frac{1}{2}}$, $H^*$-$T_b$ and $H^*$-$T_d$-Spaces}

\subsection{Definition} A topological space $X$ is said to be:
    \begin{enumerate}
        \item 
        \textbf{$H^*$-$T_1$} if for any distinct pair of points $x$ and $y$ in $X$, there exists a $H^*$-open set $U$ in $X$ containing $x$ but not $y$ and open $H^*$-open set $V$ in $X$ containing $y$ but not $x$.
        \item A \textbf{$T_\frac{1}{2}$}\cite{levine1963semi} if every $g$-closed set is closed.
        \item A \textbf{ $H^*$-$T_\frac{1}{2}$} if every $gH^*$-closed set is $H^*$-closed.
        \item A \textbf{ $H^*$-$T_b$} if every $H^*g$-closed set is closed.
        \item A \textbf{ $H^*$-$T_d$} if every $H^*g$-closed set is $g$-closed.
    \end{enumerate}

\vspace{.1cm}

   $T_1$ \hspace{.2cm} $ \Longrightarrow$ \hspace{.2cm} $ T_\frac{1}{2}$ \hspace{.2cm} $\Leftarrow$ \hspace{.2cm} $H^*$-$T_b$ \hspace{.2cm} $\Longrightarrow$ \hspace{.2cm} $H^*$-$T_d$\\

\hspace{.2cm} $\Downarrow \hspace{1.3cm} \Downarrow$ \\

\hspace{.1cm} $H^*$-$T_1$ $\Longrightarrow$ $H^*$-$T_\frac{1}{2}$
\vspace{.2cm}
\\
It can be explained by the following example:

\subsubsection{Example.} 
An \( H^* \)-\( T_b \) space need not be \( H^* \)-\( T_1 \). Let \( X = \{a, b, c, d\} \) and 
\[
\tau = \{\emptyset, X, \{a\}, \{b\}, \{a, b\}, \{a, b, c\}\}.
\]
Since \( \{a, b\} \) is not \( H^* \)-closed, the space is not \( H^* \)-\( T_1 \), and hence it is not \( T_1 \).

However, the family of all \( H^*g \)-closed sets coincides with the family of all closed sets; hence, the space is \( H^* \)-\( T_b \).

\subsection{Theorem}
\textbf{(i).} $X$ is a $T_\frac{1}{2}$ iff for each $x \in X$, $\{x\}$ is open or closed in $X$.
  \\
  \textbf{(ii).} $X$ is a $H^*$-$T_\frac{1}{2}$ space iff for each $x \in X$, $\{x\}$ is $H^*$-open or $H^*$-closed in $X$. i.e., $X$ is a $H^*$-$T_\frac{1}{2}$ iff a space $(X, \tau^{H^*})$ is $T_\frac{1}{2}$-space.

\subsection{Theorem.} 
\textbf{(i).} If $A$ is $H^*g$-closed, then $H^*$-$cl(A) - A$ does not contain non-empty closed set.\\
\textbf{(ii).} For each x$ \in X$, $\{x\}$ is closed or its complement $X - \{x\}$ is $H^*g$-closed in $X$.\\   
\textbf{(iii).} For each $x\in X$, $\{x\}$ is $H^*$-closed or its complement $X- \{x\}$ is $gH^*$-closed in $X$.

\subsection{Theorem}
\textbf{(i).} Every $H^*$-$T_b$ space is $H^*$-$T_d$ and $T_\frac{1}{2}$.\\
\textbf{(ii).} Every $T_i$ - space is $H^*$ -$T_i$, where $i $=$ 1, \frac{1}{2}$.\\
\textbf{(iii).} Every $H^*$-$T_i$ space is $H^*$-$T_\frac{1}{2}$.

\begin{proof}
    \textbf{(i).} It is obtained from \textbf{Definition 5.1}[(2), (4) and (5)] and \textbf{Remark 4.2.} (i) and (ii).\\
    \textbf{(ii).} Let $X$ be a $T_1$ (resp. $T_\frac{1}{2}$)-space and let $x\in X$. Then $\{x\}$ is closed (resp. open or closed by \textbf{Theorem 5.2.(i)}. Since every open set is $H^*$-open, $\{x\}$ is $H^*$-closed (resp. $H^*$-closed or $H^*$-open) in $X$. This implies that $X$ is $T_1$ (resp. $T_\frac{1}{2}$) by \textbf{Theorem 5.2.(ii)}. Therefore $X$ is $H^*$-$T_1$ (resp. $H^*$-$T_\frac{1}{2}$).\\
    \textbf{(iii).} Let $X$ be a $H^*$-$T_1$ space. Then $X$ is $T_1$. By \textbf{Theorem 5.3} of Levine\cite{levine1970generalized} $X$ is $T_\frac{1}{2}$ and hence $X$ is $H^*$-$T_\frac{1}{2}$.
\end{proof}

\subsection{Proposition.}
\textbf{(i).} If $X$ is $H^*$-$T_b$ then for each $x \in X$, $\{x\}$ is $H^*$-closed or open in $X$.\\
\textbf{(ii).} If $X$ is $H^*$-$T_d$ then for each $x \in X$, $\{x\}$ is $H^*$-closed or $g$-open in $X$.

\begin{proof}
   \textbf{(i).} Suppose that, for an $x \in X$, $\{x\}$ is not $H^*$-closed. By \textbf{Theorem 5.2.(iii)} and \textbf{Remark 4.2. (i) and (ii)}, $X - \{x\}$ is $H^*g$-closed set. Therefore $X - \{x\}$ closed by using assumption and hence $\{x\}$ is open.\\
   \textbf{(ii).} Suppose that, for a $x \in X$, $\{x\}$ is not closed. By \textbf{Theorem 5.3.(ii)}, $X - \{x\}$ is $H^*g$-closed set. Therefore by using assumption $X - \{x\}$ is $g$-closed and hence $\{x\}$ is $g$-open.
\end{proof}

\section{Separation Axioms $H^*$-$T_b$ and $H^*$-$T_d$ of Spaces Are Preserved Under Homomorphisms}

\subsection{Definition.} 
A map $f : X \rightarrow Y$ is said to be:
    \begin{enumerate}
        \item \textbf{pre $H^*$-closed} if for each $H^*$-closed set of $F$ of $X$, $f(F)$ is an $H^*$ closed set in $Y$.\index{pre $H^*$-closed}
        \item \textbf{$H^*$-irresolute} if for each $H^*$-closed set $F$ of $Y$, $f^{-1}(F)$ is $H^*$-closed in $X$.\index{$H^*$-irresolute}
    \end{enumerate}

\subsection{Theorem} \textbf{(i).}  A map $f : X \rightarrow Y$ pre $H^*$-closed (resp. pre $H^*$-open) iff its induced map $f : (X, \tau^{H^*})\rightarrow (Y, \sigma^{H^*})$ is a closed (resp. open) map.\\
    \textbf{(ii).} A map $f : X \rightarrow Y$ pre is $H^*$-irresolute iff its induced map \(f : (X, \tau^{H^*})\to (Y, \sigma^{H^*})\) is continuous.

\subsection{Theorem} \textbf{(i).}  A map $f : X\rightarrow Y$ is a homomorphism, then $f$ is an $H^*$-homomorphism.\\
    \textbf{(ii).} If $X$ is $H^*$-$T_b$ (resp. $H^*$-$T_d$ and $f : X\rightarrow Y$ is a homomorphism, then $Y$ is $H^*$-$T_b$ (resp. $H^*$-$T_d$).

\begin{proof}
    \textbf{(i).} Since $f : X \rightarrow Y$ is a homomorphism, then $f$ and $f ^{-1}$ are both open and $H^*$-continuous  bijection. It follows from \textbf{Theorem 4.16} of Noiri \cite{noiri1984alphacontinuous} that $f$ and $f^{-1}$ are $H^*$-irresolute. Therefore, $f$ is $H^*$-irresolute and $f$ is pre $H^*$-closed.\\
    \textbf{(ii).} Let $f : X \rightarrow Y$ is a homomorphism and let $F$ be an $H^*g$-closed set of $Y$ Then by \textbf{a.} $f^{-1} : Y \rightarrow X$ is a continuous and pre $H^*$-closed bijection. Hence by \textbf{Theorem 6.2. (i)} $f^{-1}(F)$ is $H^*g$-closed in $X$. Since $X$ is $H^*$-$T_b$ (resp. $H^*$-$T_d$), 
$f^{-1} (F)$ is closed (resp. $g$-closed) in $X$. Since $f$ is closed onto (resp. closed and continuous) map, $F$ is closed (resp. $g$-closed) by \textbf{Theorem 6.1} of Levine\cite{levine1970generalized} in $Y$. Hence $Y$ is $H^*$-$T_b$ (resp. $H^*$-$T_d$).

\end{proof}

\section*{Conclusion}

The research in topology over last two decades has reached a high level in many directions. By researching generalizations of closed sets, some new separation axioms have been founded and they turn out to be useful in the study of digital topology. Therefore, \(SC^*\) and \(H^*\)-separation axioms are defined by using \(SC^*\) and \(H^*\)-closed sets will have many possibilities of applications in digital topology and computer graphics.

%\vspace{1cm}
%\noindent {\bf Conflict of Interest Statement:} Both the authors declare that they have no conflict of interest.\\$~$\\
%\noindent {\bf Data Availability Statement:} The authors confirm that manuscript has no associated data.

\bibliographystyle{amsplain}
\bibliography{references}

\end{document}